\documentclass{article}
\usepackage{amsfonts}
\usepackage{enumitem}
\usepackage{xcolor}
\usepackage[english]{babel}

\newcommand{\keywords}[1]{\textbf{Keywords.} #1}

\usepackage{amsmath}  
\allowdisplaybreaks

\usepackage[letterpaper,top=2cm,bottom=2cm,left=3cm,right=3cm,marginparwidth=1.75cm]{geometry}

\usepackage{amsmath}
\numberwithin{equation}{section}
\usepackage{graphicx}
\usepackage[colorlinks=true, allcolors=blue]{hyperref}

\usepackage[square, numbers]{natbib}
\usepackage{amsthm}
\newtheorem{theorem}{\noindent Theorem}[section]
\newtheorem{corollary}[theorem]{\noindent Corollary}
\newtheorem{lemma}[theorem]{\noindent Lemma}
\newtheorem{proposition}[theorem]{\noindent Proposition}
\newtheorem{definition}[theorem]{\noindent Definition}
\usepackage{cleveref}
\AddToHook{env/lemma/begin}{\crefalias{theorem}{lemma}}
\AddToHook{env/corollary/begin}{\crefalias{theorem}{corollary}}

\title{On the anisotropic mean curvature flow in the Gaussian space}
\author{
{\large Di Zhang}\\
{\small School of Mathematical Sciences, Hangzhou Dianzi University, Hangzhou 310018, P.R. China}\\
{\small Email: \texttt{zhangdi5727@hdu.edu.cn}}
}
\date{}

\begin{document}
\maketitle

\begin{abstract}
In this paper, we study the anisotropic mean curvature flow in the Gaussian space endowed with a Finsler metric obtained from a Minkowski metric by a conformal deformation. We first derive the conformal transformation formulas for the fundamental tensor, the Cartan tensor, and the Chern connection of the conformally related Finsler metrics. Based on these formulas, we express the anisotropic mean curvature in the Gaussian space in terms of the corresponding anisotropic mean curvature with respect to the original Minkowski metric, and derive the associated evolution equation. Finally, we establish the short-time existence and uniqueness of smooth solutions by reformulating the geometric evolution equation as a strictly parabolic quasilinear partial differential equation.
\end{abstract}



\section{Introduction}
Let $\varphi(\cdot,t): M^n\times [0,T)\rightarrow\mathbb{R}^{n+1}$ be a family of smooth closed hypersurfaces satisfying 
\begin{align}\label{MCFeq}
    \frac{\partial}{\partial t}\varphi(x,t)=H(x,t)\mathbf{n}(x,t),
\end{align}
where $H$ denotes the mean curvature and $\mathbf{n}$ is the inward-pointing unit normal. This geometric evolution equation is known as the mean curvature flow (MCF).

The mean curvature flow was originally proposed by Mullins \cite{mul} to describe the formation of grain boundaries in annealing metals. Brakke \cite{bra} subsequently formulated the motion of submanifolds under the MCF in arbitrary codimension and constructed generalized varifold solutions existing for all time. Since then, the classical theory of the MCF for hypersurfaces has been extensively developed. In particular, Huisken \cite{hui} proved that every compact convex hypersurface in $\mathbb{R}^{n+1},\ n\geq2$, remains smooth and convex under the MCF until it contracts to a round point in finite time. He later extended this convergence result to hypersurfaces in Riemannian manifolds \cite{hui86} and studied the MCF of hypersurfaces satisfying a pinching condition in a sphere \cite{hui87}. In the higher-codimension setting, Andrews and Baker \cite{and-bak} established a convergence theorem for closed submanifolds in Euclidean space under a suitable pinching condition. Their result was subsequently extended to space forms by Baker \cite{bak} and Liu-Xu-Ye-Zhao \cite{lxyz}. In our previous work \cite{llz}, we investigated the regular geometric behavior of the MCF of submanifolds in the standard Gaussian space, which is a Riemannian manifold conformal to Euclidean space and has unbounded sectional curvature. For further developments in the theory of mean curvature flows, including local regularity, singularity formation, and problems involving phase transitions and moving interfaces, we refer the reader to [\citenum{eh, hui90, eva, ilm, and94, smo}].

Finsler geometry is a natural generalization of Riemannian geometry and exhibits a richer geometric structure. A Minkowski space is a special class of Finsler manifolds, namely, a vector space equipped with a Minkowski metric, which generalizes the Euclidean metric without requiring the metric to be quadratic. Curvature flows in Minkowski spaces have been studied in [\citenum{and00, and, xch1, xch}], where they are commonly referred to as anisotropic curvature flows. In particular, Zeng-He-Chen \cite{zhc18} introduced the concept of the mean curvature flow in the Finsler setting and provided several examples. In the setting of Minkowski spaces, they established the short-time existence and uniqueness of solutions to the mean curvature flow, proved the preservation of convexity and mean convexity, and derived comparison principles. These results provide a foundation for the further study of curvature flows in Finsler geometry and anisotropic settings. 

In this paper, we consider a Gaussian space obtained by a conformal deformation of a Minkowski space $(\mathbb{R}^{n+1},F)$. More precisely, we equip $\mathbb{R}^{n+1}$ with the Finsler metric
$$\widetilde F(x,y)=e^{-\frac{|x|^2}{n}}F(y),$$
where $|x|$ denotes the standard Euclidean norm of $x\in\mathbb{R}^{n+1}$. The resulting Finsler manifold $(\mathbb{R}^{n+1},\widetilde F)$ is called the Gaussian space associated with $(\mathbb{R}^{n+1},F)$. Equivalently, we have $\widetilde F=e^{\sigma(x)}F$, where $\sigma(x)=-\frac{|x|^2}{n}$. 

Motivated by the work of Zeng-He-Chen \cite{zhc18}, as well as by our previous study of the mean curvature flow of submanifolds in Gaussian space in the Riemannian setting \cite{llz}, we investigate the anisotropic mean curvature flow in the Gaussian space defined above. More precisely, we consider the evolution equation
\begin{align}
    \begin{cases}
        \partial_t\varphi(p,t)
        =\widetilde{H}(p,t)\widetilde{\mathbf{n}}_0(p,t),\\
        \varphi(p,0)=\varphi_0(p),
    \end{cases}
\end{align}
where $\widetilde{H}(p,t)$ and $\widetilde{\mathbf{n}}_0(p,t)$ denote, respectively, the anisotropic mean curvature and the inward-pointing unit normal of the hypersurface $\varphi_t$ with respect to the Finsler metric $\widetilde{F}$ at $p\in M$.

Following the approach used in our earlier work on the conformal mean curvature flow in the Riemannian setting (see the introduction of \cite{lz}), we study the anisotropic mean curvature flow in the Gaussian space by expressing its geometric quantities in terms of those associated with the underlying Minkowski metric $F$. To this end, we derive conformal transformation formulas for the fundamental tensor, the Cartan tensor, and the Chern connection. These formulas yield an explicit expression for the anisotropic mean curvature in the Gaussian Finsler space in terms of the corresponding quantities associated with $F$, and enable us to reformulate the flow in terms of the anisotropic mean curvature $H$ and the unit normal $\mathbf{n}$ with respect to $F$, rather than working directly with $\widetilde{H}$ and $\widetilde{\mathbf{n}}$. This gives the following initial value problem:
\begin{align}
    \begin{cases}
        \partial_t\varphi(p,t)
        =e^{\frac{2}{n}\left|\varphi\right|^2}\left(H+2\langle\varphi,\mathbf{n}\rangle+\frac{1}{n}\widehat{I}(\nabla^\mathbf{n}|\varphi|^2)\right)\mathbf{n},\\
        \varphi(p,0)=\varphi_0(p),
    \end{cases}
\end{align}
where $\langle\varphi,\mathbf{n}\rangle=\sum_{\kappa=1}^{n+1}\varphi^\kappa\mathbf{n}^\kappa$ is the standard Euclidean inner product of $\varphi$ and $\mathbf{n}$, and
\begin{align*}
    \widehat I(\nabla^\mathbf{n}|\varphi|^2)=\widehat g^{ij}C_{\mathbf{n}}\left(\varphi_*(\partial_i),\varphi_*(\partial_j),\nabla^\mathbf{n}|\varphi|^2\right)=2\widehat{g}^{ij}C_{\beta\gamma}^{\kappa}(\mathbf{n})\varphi^\kappa\varphi^\beta_i\varphi^\gamma_j
\end{align*}
is the mean Cartan torsion evaluated at $\nabla^\mathbf{n}|\varphi|^2$. We emphasize that $\mathbf{n}$ is the unit normal with respect to the Minkowski metric $F$, and hence its Euclidean length is not necessarily equal to one. In particular, the term $\langle\varphi,\mathbf{n}\rangle$ represents one half of the Euclidean directional derivative $\mathbf{n}(|\varphi|^2)$.

Using this formulation, we also establish short-time existence and uniqueness (see Theorem \ref{mainthm}), up to the reparametrization specified in Corollary \ref{coro+1}, for the flow starting from any smooth compact initial hypersurface. This is done by reducing the normal-velocity equation to a uniformly strictly parabolic quasilinear equation for the normal graph function.

The remainder of this paper is organized as follows. In Section 2, we recall the fundamental concepts and formulas from Finsler geometry needed in the sequel. In Section 3, we derive the conformal transformation formulas relating the Finsler metrics $F$ and $\widetilde{F}$. In Section 4, we obtain the anisotropic mean curvature flow equation in Gaussian space. Finally, in Section 5, we reformulate the geometric evolution equation as a strictly parabolic second-order quasilinear partial differential equation and establish the short-time existence and uniqueness of smooth solutions.

\section{Preliminaries}
\subsection{Finsler manifolds}
Throughout this paper, unless otherwise stated, we use the following convention for index ranges:
$$1\le i,j,k,\cdots\le n, \ 1\le\alpha,\beta,\gamma,\cdots\le n+1, \ 1\le I,J,\cdots\le 2n.$$
Let $M$ be an $n$-dimensional oriented smooth manifold. We denote by $TM$ its tangent bundle and by $(x,y)=(x^1,\cdots x^n,y^1,\cdots y^n)$ the induced local coordinates on $TM$. A \textit{Finsler metric} on $M$ is a function $F:TM\rightarrow [0,+\infty)$ satisfying the following conditions:
\begin{enumerate}[label=(\roman*)]
    \item $F$ is smooth on the slit tangent bundle $TM\backslash{\{0\}}$;
    \item $F(x,\lambda y)=\lambda F(x,y)$ for all $\lambda >0$;
    \item The $n\times n$ matrix $(g_{ij}):=\left(\frac{1}{2}\frac{\partial^2F^2}{\partial y^i\partial y^j}\right)$ is positive definite at every point of $TM\backslash{\{0\}}$.
\end{enumerate}
The pair $(M,F)$ is called a \textit{Finsler manifold}.

The \textit{fundamental tensor} $g$ and the \textit{Cartan tensor} $C$ of $F$ are defined by (see \cite{book-xia})
$$g=g_{ij}(x,y)\mathrm{d}x^i\otimes\mathrm{d}x^j,\quad C=C_{ijk}(x,y)\mathrm{d}x^i\otimes\mathrm{d}x^j\otimes\mathrm{d}x^k,$$
where
$$g_{ij}(x,y):=\frac{1}{2}[F^2(x,y)]_{y^iy^j}=FF_{y^iy^j}+F_{y^i}F_{y^j},\quad\text{and}\quad C_{ijk}:=\frac{1}{4}[F^2(x,y)]_{y^iy^jy^k}=\frac{1}{2}\frac{\partial{g_{ij}}}{\partial{y^k}}.$$

By Euler's theorem on homogeneous functions, we have
\begin{align*}
    F=F_{y^i}y^i,\quad F_{y^iy^j}y^i=0,\quad FF_{y^i}=g_{ij}(x,y)y^j,
\end{align*}
and
\begin{align*}
    F^2(x,y)=g_{ij}(x,y)y^iy^j,\quad C_{ijk}(x,y)y^i=0.
\end{align*}
In particular, $F$ is Riemannian if and only if its Cartan tensor vanishes identically, that is, $C_{ijk}=0$ for all $1\le i,j,k\le n$.

The \textit{mean Cartan torsion} $I$ is defined by $I(y):=\mathrm{tr}_{g}C(y,\cdot,\cdot)$ for $y\in T_xM$. In local coordinates,
\begin{align*}
    I=I_i(x,y)\mathrm{d}x^i,\quad I_i=g^{jk}C_{ijk}.
\end{align*}
It follows that $I_iy^i=0$. Moreover,
$$I_i=\frac{1}{2}g^{jk}\frac{\partial g_{ij}}{\partial y^k}=\frac{1}{2}g^{jk}\frac{\partial g_{jk}}{\partial y^i}=\frac{\partial}{\partial y^i}\left(\mathrm{log}\sqrt{\mathrm{det}(g_{ij})}\right).$$

The geodesics of $F$ are locally characterized by the system of second-order ordinary differential equations
$$\frac{\mathrm{d}^2x^i}{\mathrm{d}t^2}+G^i\left(x,\frac{\mathrm{d}x}{\mathrm{d}t}\right)=0,$$
where 
$$G^i(x,y)=\frac{1}{2}g^{il}\left(\frac{\partial{g_{jl}}}{\partial{x^k}}+\frac{\partial{g_{kl}}}{\partial{x^j}}-\frac{\partial{g_{jk}}}{\partial{x^l}}\right)y^jy^k.$$
The functions $G^i$ are called the \textit{geodesic coefficients} of $F$. The coefficients of the nonlinear connection are given by $N^i_j=\frac{1}{2}\frac{\partial{G^i}}{\partial{y^j}}$ (see \cite{bcs}).

Given a Finsler metric $F$, its {\textit{dual metric}} $F^*:T^*M\rightarrow[0,+\infty)$ is defined by
$$F^*(x,\xi):=\sup_{y\in T_xM\backslash{\{0\}}}\frac{\xi(y)}{F(x,y)}.$$
The dual metric $F^*$ is also a Finsler metric on $M$.

The canonical projection $\pi:TM\rightarrow M$ induces the pull-back bundles $\pi^*TM$ and $\pi^*T^*M$ over $TM\backslash{\{0\}}$. The \textit{distinguished section} $\ell$ of $\pi^*TM$ is defined by (see \cite{bcs})
$$\ell=\ell_{(x,y)}:=\frac{y^i}{F(y)}\frac{\partial}{\partial x^i}=\frac{y^i}{F}\frac{\partial}{\partial x^i}=:\ell^i\frac{\partial}{\partial x^i}.$$
Its natural dual is the \textit{Hilbert form} $\omega$, a section of $\pi^*T^*M$, given by
$$\omega=\omega_{(x,y)}:=F_{y^i}(x,y)\mathrm{d}x^i=F_{y^i}\mathrm{d}x^i.$$
Both $\ell$ and $\omega$ are globally defined on $TM\backslash{\{0\}}$.

It is well known that there exists a unique \textit{Chern connection} $\nabla$ on $\pi^*TM$, characterized by the following two properties (see \cite{bcs}). Writing
$\nabla\frac{\partial} {\partial{x^i}}=\omega_i^j\frac{\partial}{\partial{x^j}}$, its connection forms satisfy
$$\omega_i^j\wedge\mathrm{d}x^i=0,$$
and
$$\mathrm{d}g_{ij}-g_{kj}\omega_i^k-g_{ik}\omega_j^k=2A_{ijk}\frac{\delta y^k}{F},$$
where
$$\delta y^k:=\mathrm{d}y^k+N^k_j\mathrm{d}x^j=\mathrm{d}y^k+y^j\omega_j^k,\quad A_{ijk}=FC_{ijk}.$$

For a vector field $X=X^i\frac{\partial}{\partial x^i}$, its \textit{covariant derivative} along a vector $v=v^i\frac{\partial}{\partial x^i}$, with reference vector $w\in T_xM\backslash{\{0\}}$, is defined by
\begin{align*}
    D^w_vX(x):=\left\{v^j\frac{\partial X^i}{\partial x^j}(x)+\Gamma_{jk}^i(w)v^jX^k(x)\right\}\frac{\partial}{\partial x^i},
\end{align*}
where $\Gamma_{jk}^i$ are the connection coefficients of the Chern connection.

The \textit{Legendre transformation} $\mathcal{L}:TM\rightarrow T^*M$ is defined by
\begin{align*}
    \mathcal{L}(x,y):=
    \begin{cases}
        g_y(y,\cdot),\ \ \ \ (x,y)\in TM\backslash{\{0\}},\\
        0, \quad\quad\quad\ \ \ (x,y)=0\in TM,
    \end{cases}
\end{align*}
where $g_y(y,\cdot):=g_{ij}(x,y)y^i\mathrm{d}x^j$. For a smooth function $f:M\rightarrow\mathbb{R}$, its \textit{Finsler gradient} is given by $\nabla{f}(x):=\mathcal{L}^{-1}(\mathrm{d}f(x))$. In local coordinates, this becomes
\begin{align*}
    \nabla{f}(x):=
    \begin{cases}
        g^{ij}(x,\nabla{f})\frac{\partial{f}}{\partial{x^j}}\frac{\partial}{\partial{x^i}},\quad\mathrm{d}f(x)\neq 0,\\
        0, \qquad\qquad\qquad\qquad\mathrm{d}f(x)=0.
    \end{cases}
\end{align*}
We denote by $M_f:=\left\{x\in M\mid\mathrm{d}f(x)\neq0\right\}$ the regular set of $f$. For $x\in M_f$, the Hessian of $f$ is defined by 
$$\nabla^2f(v):=D^{\nabla f}_v\nabla f(x),\quad v\in T_xM,$$
where the covariant derivative is taken with reference vector $\nabla f$ (see \cite{yhs}).

For a nonzero reference vector $v\in T_xM$, we also define the \textit{reference-vector gradient} $\nabla^vf(x)$ of $f$ with respect to the inner product $g_v$. More precisely, $\nabla^vf(x)$ is the unique vector in $T_xM$ satisfying
\begin{align*}
    g_v(\nabla^vf(x),X)=\mathrm{d}f(x)\quad\text{for every}\quad X\in T_xM.
\end{align*}
In local coordinates,
\begin{align*}
    (\nabla^vf)^i=g^{ij}(x,v)\frac{\partial f}{\partial x^j}.
\end{align*}
Here $g^{ij}(x,v)$ denotes the inverse matrix of $g_{ij}(x,v)$. This reference-vector gradient is distinct from the canonical Finsler gradient $\nabla f=\mathcal{L}^{-1}(\mathrm{d}f)$. Whenever $\mathrm{d}f\neq0$, they are related by $\nabla f=\nabla^{\nabla f}f$.

\subsection{Anisotropic hypersurfaces}\label{Sec2.2}
Let $\varphi:M^n\rightarrow (N^{n+1}, F)$ be an embedded hypersurface of the Finsler manifold $(N,F)$, where $M$ is an $n$-dimensional oriented smooth manifold. Let $(x^\alpha, y^\alpha)$ and $(u^i)$ denote local coordinates on $TN$ and $M$, respectively. Following \cite{zhc18}, we define the \textit{conormal bundle} of the hypersurface by
$$\mathcal{V}(M)=\left\{\xi\in\varphi^{-1}T^*N\mid\xi(\mathrm{d}\varphi(X))=0,\quad\forall X\in\Gamma(TM)\right\}.$$ 
In the possibly nonreversible case, the two conormal rays must be normalized separately. Let $\nu\in\mathcal{V}(M)$ be a smooth conormal field satisfying $F^*(\nu)=1$. We define $\nu^+=\nu$ and $\nu^-=\frac{-\nu}{F^*(-\nu)}$. Then $F^*(\nu^+)=F^*(\nu^-)=1$. The corresponding $F$-unit normal vector fields are given by $\mathbf{n}^+=\mathcal{L}^{-1}(\nu^+)$ and $\mathbf{n}^-=\mathcal{L}^{-1}(\nu^-)$. In general, $\mathbf{n}^-\neq-\mathbf{n}^+$ unless $F$ is reversible. In what follows, we fix one coorientation and denote the corresponding $F$-unit normal field simply by $\mathbf{n}$. 

With this choice of coorientation, let $\mathbf{n}=\mathcal{L}^{-1}(\nu)$ be the unit normal field along $M$, and set $\widehat{g}:=\varphi^*g_{\mathbf{n}}$. Then $(M,\widehat{g})$ is a Riemannian manifold. We refer to $(M,\widehat g)$ as the \textit{anisotropic hypersurface} induced by $\varphi$, in order to distinguish its induced Riemannian metric from the ambient Finsler metric. In local coordinates, 
$\widehat g_{ij}=g_{\alpha\beta}(\mathbf{n})\varphi^\alpha_i\varphi^\beta_j$, where
$$x^\alpha=\varphi^\alpha(u),\quad y=\mathbf{n},\quad \varphi^\alpha_i=\frac{\partial\varphi^\alpha}{\partial u^i}.$$

The following \textit{Gauss-Weingarten formulas} hold, analogously to the Riemannian case:
\begin{align}
    D^\mathbf{n}_X\varphi_*(Y)=\varphi_*(\widehat\nabla_XY)+h(X,Y)\mathbf{n},\quad D^\mathbf{n}_X\mathbf{n}=-\varphi_*(A(X)),
\end{align}
for $X,Y\in\Gamma(TM)$, where
\begin{align}
    h(X,Y):=\widehat g(AX,Y)=g_\mathbf{n}\big(\mathbf{n},D^\mathbf{n}_X\varphi_*(Y)\big).
\end{align}
Here $A:T_uM\rightarrow T_uM$ is a linear map that is self-adjoint with respect to $\widehat g$, and is called the \textit{Weingarten transformation}(or \textit{shape operator}) of the anisotropic hypersurface. 

The induced connection $\widehat\nabla$ is a torsion-free linear connection on $M$ satisfying
\begin{align}
    (\widehat\nabla_X\widehat g)(Y,Z)=2C_\mathbf{n}(D^\mathbf{n}_X\mathbf{n},Y,Z)=-2\widehat{C}(A(X),Y,Z),
\end{align}
where $\widehat{C}(A(X),Y,Z)=C_\mathbf{n}\left(\varphi_*(A(X)),\varphi_*(Y),\varphi_*(Z)\right)$. Consequently, $\widehat\nabla$ is generally not the Levi-Civita connection of $(M,\widehat g)$, unless the Cartan tensor vanishes along the relevant directions.

The bilinear form $h$ is symmetric, that is, $h(X,Y)=h(Y,X)$, and is called the \textit{anisotropic second fundamental form} of $M$ in $(N,F)$. The eigenvalues of $A$, denoted by $k_1,k_2,\cdots,k_n,$ are called the \textit{anisotropic principal curvatures} of $M$ with respect to the chosen unit normal $\mathbf{n}$. If there exists a function $k$ such that $A=k\mathrm{Id}$, or equivalently $k_1=k_2=\cdots=k_n$, then $M$ is called an \textit{anisotropic totally umbilic} hypersurface of $(N,F)$.

Following [\citenum{hys},\citenum{shen}], we define the \textit{anisotropic mean curvature} of $M$ in $(N,F)$ by
$$H=\mathrm{tr}_{\widehat g}h=\sum_{i=1}^nk_i,$$
without dividing by $n$. When $(N,F)$ is a Minkowski space, this definition agrees with the anisotropic mean curvature used, for example, in [\citenum{and},\citenum{xch}].

\section{Conformal transformation formulas for Finsler metrics}
In this section, we derive several conformal transformation formulas for two conformally related Finsler metrics $F$ and $\widetilde F=e^{\sigma(x)}F$ on an $(n+1)$-dimensional smooth manifold $N$, where $\sigma$ is a smooth function on $N$. For the definition of conformally related Finsler metrics, we refer the reader to the introduction of \cite{ccz}. Most of the results presented in this section are adapted from \cite{has76}. For the convenience of the reader, we restate these resutlts using our notation and provide detailed proofs where necessary, thereby making the exposition self-contained.

\begin{lemma}{\textup{(\cite{has76})}}\label{lem1}
    Let $F$ and $\widetilde F=e^{\sigma}F$ be two conformally related Finsler metrics on $N$. Then
    \begin{align*}
        \widetilde g_{\alpha\beta}=e^{2\sigma}g_{\alpha\beta},
        \quad\widetilde g^{\alpha\beta}=e^{-2\sigma}g^{\alpha\beta},
    \end{align*}
    and
    \begin{align*}
        \widetilde\ell^\alpha=e^{-\sigma}\ell^\alpha,\quad \widetilde\ell_\alpha=e^{\sigma}\ell_\alpha,\quad 
        \widetilde y_\alpha=e^{2\sigma}y_\alpha,
    \end{align*}
    where $(g^{\alpha\beta})$ is the inverse matrix of $(g_{\alpha\beta})$, and
    \begin{align*}
        \ell_\alpha:=F_{y^\alpha}=g_{\alpha\beta}\ell^\beta=\frac{y_\alpha}{F},\quad y_\alpha:=g_{\alpha\beta}y^\beta.   
    \end{align*}
\end{lemma}

\begin{lemma}{\textup{(\cite{has76})}}\label{lem2}
    Under the same conformal change, the Cartan tensor satisfies
    \begin{align}
        \widetilde C_{\alpha\beta\gamma}=e^{2\sigma}C_{\alpha\beta\gamma},
        \quad
        \widetilde C^{\alpha}_{\beta\gamma}=C^{\alpha}_{\beta\gamma},
        \quad
        \widetilde C^{\alpha\beta}_{\gamma}=e^{-2\sigma}C^{\alpha\beta}_{\gamma},
    \end{align}
    where 
    \begin{align*}
        C_{\alpha\beta\gamma}:=\frac{1}{2}\frac{\partial g_{\alpha\beta}}{\partial y^\gamma},\quad C^{\alpha}_{\beta\gamma}:=g^{\alpha\kappa}C_{\kappa\beta\gamma},\quad C^{\alpha\beta}_{\gamma}:=g^{\beta\eta}C^{\alpha}_{\eta\gamma}=-\frac{1}{2}\frac{\partial g^{\alpha\beta}}{\partial y^\gamma}. 
    \end{align*}
\end{lemma}

We next introduce a conformally invariant tensor that will be used in the subsequent computations.
\begin{lemma}{\textup{(\cite{has76})}}\label{lem3}
    The tensor $B_{\alpha\beta}:=(\log F^2)_{y^\alpha y^\beta}$ is conformally invariant, that is, 
    $\widetilde B_{\alpha\beta}=B_{\alpha\beta}$. Moreover,
    \begin{align}\label{bab}
        B_{\alpha\beta}=\frac{2}{F^2}(g_{\alpha\beta}-2\ell_\alpha\ell_\beta).
    \end{align}
    The inverse matrix $(B_{\alpha\beta})$ is given by $(B^{\alpha\beta})$ given by
    \begin{align}\label{bab-1}
        B^{\alpha\beta}=\frac{F^2}{2}(g^{\alpha\beta}-2\ell^\alpha\ell^\beta),
    \end{align}
    and is also conformally invariant.  
\end{lemma}

\begin{proof}
    Since $\widetilde F^2=e^{2\sigma}F^2$, we have
    $$\log\widetilde F^2=2\sigma+\log F^2.$$
    Differentiating twice with respect to $y^\alpha$ and $y^\beta$, and using the fact that $\sigma$ is independent of $y$, yields
    $$(\log\widetilde F^2)_{y^\alpha y^\beta}=(\log F^2)_{y^\alpha y^\beta}.$$
    Hence $\widetilde B_{\alpha\beta}=B_{\alpha\beta}$.

    To obtain the explicit expression, note that
    $$B_{\alpha\beta}=2(\log F)_{y^\alpha y^\beta}.$$
    Using $\ell_\alpha:=F_{y^\alpha}$, we find
    $$(\log F)_{y^\alpha y^\beta}=\frac{1}{F}\frac{\partial\ell_\alpha}{\partial y^\beta}-\frac{1}{F^2}\ell_\alpha\ell_\beta.$$
    On the other hand, $g_{\alpha\beta}=F\frac{\partial\ell_\alpha}{\partial y^\beta}+\ell_\alpha\ell_\beta$, and therefore
    $$\frac{\partial\ell_\alpha}{\partial y^\beta}=\frac{g_{\alpha\beta}-\ell_\alpha\ell_\beta}{F}.$$
    Substitution gives
    $$B_{\alpha\beta}=2(\log F)_{y^\alpha y^\beta}=\frac{2}{F^2}(g_{\alpha\beta}-2\ell_\alpha\ell_\beta),$$
    which proves \eqref{bab}.

    It remains to verify that \eqref{bab-1} defines the inverse of $B_{\alpha\beta}$. Using $g_{\alpha\beta}\ell^\beta=\ell_\alpha$ and $\ell_\beta\ell^\beta=1$, we obtain
    \begin{align*}
        B_{\alpha\beta}B^{\beta\gamma}
        =&(g_{\alpha\beta}-2\ell_\alpha\ell_\beta)(g^{\beta\gamma}-2\ell^\beta\ell^\gamma)\\
        =&\delta_\alpha^\gamma-2\ell_\alpha\ell^\gamma-2\ell_\alpha\ell^\gamma+4\ell_\alpha\ell^\gamma\\
        =&\delta_\alpha^\gamma.
    \end{align*}
    Thus $(B^{\alpha\beta})$ is the inverse of $(B_{\alpha\beta})$. Its conformal invariance follows from that of $B_{\alpha\beta}$.
\end{proof}

\begin{lemma}{\textup{(\cite{has76})}}\label{lem4}
    The formal Christoffel symbols of the fundamental tensor are defined by
    \begin{align}\label{gamma}
        \gamma^\alpha_{\beta\gamma}:=\frac{1}{2}g^{\alpha\eta}\left(\frac{\partial g_{\beta\eta}}{\partial x^\gamma}+\frac{\partial g_{\gamma\eta}}{\partial x^\beta}-\frac{\partial g_{\beta\gamma}}{\partial x^\eta}\right).
    \end{align}
    Under the conformal change $\widetilde{g}_{\alpha\beta}=e^{2\sigma}g_{\alpha\beta}$, they transform according to
    \begin{align}\label{gamma1}
        \widetilde\gamma^\alpha_{\beta\gamma}=\gamma^\alpha_{\beta\gamma}-\left(g_{\beta\gamma}\sigma^\alpha-\delta^\alpha_\beta\sigma_\gamma-\delta^\alpha_\gamma\sigma_\beta\right),
    \end{align}
    where
    $$\sigma_\kappa:=\frac{\partial\sigma}{\partial x^\kappa},\quad \sigma^\alpha:=g^{\alpha\kappa}\sigma_\kappa,$$ and $\delta^\alpha_\beta$ denotes the Kronecker delta.
\end{lemma}

\begin{proof}
Substituting $\widetilde{g}_{\alpha\beta}=e^{2\sigma}g_{\alpha\beta}$ and $\widetilde{g}^{\alpha\eta}=e^{-2\sigma}g^{\alpha\eta}$ into \eqref{gamma}, and using
$$\frac{\partial\widetilde{g}_{\beta\eta}}{\partial x^\gamma}=e^{2\sigma}(\frac{\partial g_{\beta\eta}}{\partial x^{\gamma}}+2g_{\beta\eta}\sigma_\gamma),$$
we obtain the stated transformation formula after simplification.
\end{proof}

\begin{lemma}\label{lem5}
    Set $G^\alpha:=\gamma^\alpha_{\beta\gamma}y^\beta y^\gamma$, $B^\alpha:=B^{\alpha\kappa}\sigma_\kappa$. Then
    \begin{align}\label{G}
        \widetilde G^\alpha=G^\alpha-2B^\alpha.
    \end{align} 
\end{lemma}

\begin{proof}
    By definition,
    $$\widetilde{G}=\widetilde{\gamma}_{\beta\gamma}^\alpha y^\beta y^\gamma.$$
    Using \eqref{gamma1}, we obtain
    \begin{align*}
        \widetilde{G}^\alpha=&(\gamma_{\beta\gamma}^{\alpha}-g_{\beta\gamma}\sigma^\alpha+\delta_\beta^\alpha\sigma_\gamma+\delta_\gamma^\alpha\sigma_\beta)y^\beta y^\gamma\\
        =&G^\alpha-g_{\beta\gamma}y^\beta y^\gamma\sigma^\alpha+\sigma_\gamma y^\gamma y^\alpha+\sigma_\beta y^\beta y^\alpha\\
        =&G^\alpha-F^2\sigma^\alpha+2\sigma_0y^\alpha,
    \end{align*}
    where $\sigma_0:=\sigma_\gamma y^\gamma$.

    By \eqref{bab-1},
    \begin{align*}
        B^{\alpha}=B^{\alpha\kappa}\sigma_\kappa=\frac{F^2}{2}(g^{\alpha\kappa}-2\ell^\alpha\ell^\kappa)\sigma_\kappa
        =\frac{F^2}{2}\sigma^\alpha-F^2\ell^\alpha\ell^\kappa\sigma_\kappa=\frac{F^2}{2}\sigma^\alpha-\sigma_0y^\alpha.
    \end{align*}
    Consequently,
    $$2B^\alpha=F^2\sigma^\alpha-2\sigma_0y^\alpha,$$
    and hence
    $$\widetilde{G}^\alpha=G^\alpha-2B^\alpha.$$
    This proves \eqref{G}.   
\end{proof}

\begin{lemma}{\textup{(\cite{has76})}}\label{lem6}
   A direct computation yields the identity
   \begin{align}\label{Bar}
       \frac{\partial B^{\alpha\gamma}}{\partial y^\beta}=y_\beta g^{\alpha\gamma}-\delta_\beta^\alpha y^\gamma-\delta_\beta^\gamma y^\alpha-F^2C_\beta^{\alpha\gamma}.
   \end{align}
\end{lemma}

\begin{proof}
    Recall from \eqref{bab-1} that
    $$B^{\alpha\gamma}=\frac{F^2}{2}(g^{\alpha\gamma}-2\ell^\alpha\ell^\gamma).$$
    Differentiating with respect to $y^\beta$, we obtain
    $$\frac{\partial B^{\alpha\gamma}}{\partial y^\beta}=\frac{1}{2}\frac{\partial F^2}{\partial y^\beta}(g^{\alpha\gamma}-2\ell^\alpha\ell^\gamma)+\frac{F^2}{2}(\frac{\partial g^{\alpha\gamma}}{\partial y^\beta}-2\frac{\partial\ell^\alpha}{\partial y^\beta}\ell^\gamma-2\ell^\alpha\frac{\partial\ell^\gamma}{\partial y^\beta}).$$
    We evaluate the terms on the right-hand side separately. Since $F^2=g_{\mu\nu}y^\mu y^\nu$, we have
    $$\frac{\partial F^2}{\partial y^\beta}=2g_{\mu\beta}y^\mu=2y_\beta.$$
    Moreover, by the definition $\ell^\alpha=y^\alpha/F$,
    $$\frac{\partial\ell^\alpha}{\partial y^\beta}=\frac{\delta_\beta^\alpha}{F}-\frac{y^\alpha}{F^2}\frac{\partial F}{\partial y^\beta}=\frac{\delta_\beta^\alpha}{F}-\frac{y^\alpha\ell_\beta}{F^2}=\frac{1}{F}(\delta_\beta^\alpha-\ell^\alpha\ell_\beta).$$
    Finally, recalling that $C_\beta^{\alpha\gamma}:=-\frac{1}{2}\frac{\partial g^{\alpha\gamma}}{\partial y^\beta}$, we obtain
    $$\frac{\partial g^{\alpha\gamma}}{\partial y^\beta}=-2C_\beta^{\alpha\gamma}.$$
    Substituting these identities into the expression for $\frac{\partial B^{\alpha\gamma}}{\partial y^\beta}$ gives
    \begin{align*}
        \frac{\partial B^{\alpha\gamma}}{\partial y^\beta}=&y_\beta(g^{\alpha\gamma}-2\ell^\alpha\ell^\gamma)
        -F^2C_\beta^{\alpha\gamma}-F(\delta_\beta^\alpha-\ell^\alpha\ell_\beta)\ell^\gamma-F\ell^\alpha(\delta_\beta^\gamma-\ell^\gamma\ell_\beta)\\
        =&y_\beta g^{\alpha\gamma}-2y_\beta\ell^\alpha\ell^\gamma-F^2C_\beta^{\alpha\gamma}-F\delta_\beta^\alpha\ell^\gamma-F\ell^\alpha\delta_\beta^\gamma+2F\ell^\alpha\ell_\beta\ell^\gamma.
    \end{align*}
    Using $F\ell^\alpha=y^\alpha$ and $F\ell_\beta=y_\beta$, the terms involving $\ell^\alpha$ combine as
    $$-2y_\beta\ell^\alpha\ell^\gamma+2F\ell^\alpha\ell_\beta\ell^\gamma=-2y_\beta\ell^\alpha\ell^\gamma+2y_\beta\ell^\alpha\ell^\gamma=0.$$
    Hence,
    \begin{align*}
        \frac{\partial B^{\alpha\gamma}}{\partial y^\beta}
        =&y_\beta g^{\alpha\gamma}-\delta_\beta^\alpha y^\gamma-\delta_\beta^\gamma y^\alpha-F^2C_\beta^{\alpha\gamma},
    \end{align*}
    which proves \eqref{Bar}.
\end{proof}

\begin{lemma}{\textup{(\cite{has76})}}\label{lem7}
    The following identity holds:
    \begin{align}\label{Bab}
        B^\alpha_\beta:=\frac{\partial B^\alpha}{\partial y^\beta}=y_\beta\sigma^\alpha-\delta^\alpha_\beta\sigma_0-y^\alpha\sigma_\beta-F^2C_\beta^\alpha,
    \end{align}
    where $\sigma_0:=\sigma_\gamma y^\gamma$ and $C_\beta^\alpha:=C_\beta^{\alpha\gamma}\sigma_\gamma$.
\end{lemma}

\begin{proof}
    Recall that $B^\alpha=B^{\alpha\kappa}\sigma_\kappa$, where $\sigma_\kappa=\sigma_\kappa(x)$ is independent of $y$. Differentiating with respect to $y^\beta$, we obtain
    $$B^\alpha_\beta=\frac{\partial B^\alpha}{\partial y^\beta}=\frac{\partial B^{\alpha\kappa}}{\partial y^\beta}\sigma_\kappa.$$
    By Lemma \ref{lem6}, this becomes
    \begin{align*}
        B_\beta^\alpha=&(y_\beta g^{\alpha\kappa}-\delta_\beta^\alpha y^\kappa-\delta_\beta^\kappa y^\alpha-F^2C_\beta^{\alpha\kappa})\sigma_\kappa\\
        =&y_\beta(g^{\alpha\kappa}\sigma_\kappa)-\delta_\beta^\alpha(y^\kappa\sigma_\kappa)-y^\alpha(\delta_\beta^\kappa\sigma_\kappa)-F^2(C_\beta^{\alpha\kappa}\sigma_\kappa).
    \end{align*}
    Using
    $$g^{\alpha\kappa}\sigma_\kappa=\sigma^\alpha,\quad y^\kappa\sigma_\kappa=\sigma_0,\quad \delta_\beta^\kappa\sigma_\kappa=\sigma_\beta,\quad C_\beta^{\alpha\kappa}\sigma_\kappa=C_\beta^\alpha,$$
    we conclude that
    $$B_\beta^\alpha=y_\beta\sigma^\alpha-\delta_\beta^\alpha\sigma_0-y^\alpha\sigma_\beta-F^2C_\beta^\alpha.$$
    This proves \eqref{Bab}.
\end{proof}

\begin{lemma}\label{lem8}
    Differentiating \eqref{G} with respect to $y^\beta$ yields
    \begin{align}\label{N}
        \widetilde N^\alpha_\beta=N^\alpha_\beta-B^\alpha_\beta,
    \end{align}
    where $N^\alpha_\beta:=\frac{\partial}{\partial y^\beta}\left(\frac{1}{2}\gamma^\alpha_{\lambda\mu}y^\lambda y^\mu\right)$.
\end{lemma}

\begin{proof}
    By definition, $N^\alpha_\beta=\frac{1}{2}\frac{\partial G^\alpha}{\partial y^\beta}$. Differentiating \eqref{G} with respect to $y^\beta$, and using $B_\beta^\alpha=\frac{\partial B^\alpha}{\partial y^\beta}$ from \eqref{Bab}, we obtain
    \begin{align*}
        \widetilde N^\alpha_\beta=\frac{1}{2}\frac{\partial\widetilde G^\alpha}{\partial y^\beta}=\frac{1}{2}\frac{\partial G^\alpha}{\partial y^\beta}-\frac{\partial B^\alpha}{\partial y^\beta}=N^\alpha_\beta-B^\alpha_\beta.
    \end{align*}
    This proves \eqref{N}.
\end{proof}

We now derive the conformal transformation formulas for the Chern connection coefficients.
\begin{lemma}\label{lem9}
    Recall that the Chern connection coefficients are given by
    \begin{align}\label{Gamma1}
        \Gamma_{\beta\gamma}^\alpha=\gamma_{\beta\gamma}^\alpha-g^{\alpha\eta}\left(C_{\eta\beta\delta}N_\gamma^\delta-C_{\beta\gamma\delta}N_\eta^\delta+C_{\gamma\eta\delta}N_\beta^\delta\right).
    \end{align}
    Under the conformal change $\widetilde{g}_{\alpha\beta}=e^{2\sigma}g_{\alpha\beta}$, these coefficients transform as
    \begin{align}\label{Gamma2}
        \widetilde\Gamma_{\beta\gamma}^\alpha=\Gamma_{\beta\gamma}^\alpha-\left(g_{\beta\gamma}g^{\alpha\kappa}\sigma_{\kappa}-\delta_{\beta}^{\alpha}\sigma_{\gamma}-\delta_{\gamma}^{\alpha}\sigma_{\beta}\right)
        +g^{\alpha\eta}\left(C_{\eta\beta\delta}B_{\gamma}^{\delta}-C_{\beta\gamma\delta}B_{\eta}^{\delta}+C_{\gamma\eta\delta}B_{\beta}^{\delta}\right),
    \end{align}
    where $B_{\beta}^{\alpha}$ is given by 
    \begin{align}\label{B}
        B_{\beta}^{\alpha}(y)=g_{\beta\gamma}(y)y^{\gamma}g^{\alpha\kappa}(y)\sigma_{\kappa}-\delta_{\beta}^{\alpha}\sigma_{\kappa}y^{\kappa}-y^{\alpha}\sigma_{\beta}-F^2g^{\alpha\lambda}(y)g^{\kappa\mu}(y)C_{\lambda\mu\beta}\sigma_{\kappa}.
    \end{align}
\end{lemma}

\begin{proof}
    By the definition \eqref{Gamma1}, we have
    $$\widetilde\Gamma_{\beta\gamma}^\alpha=\widetilde\gamma_{\beta\gamma}^\alpha-\widetilde{g}^{\alpha\eta}\left(\widetilde{C}_{\eta\beta\delta}\widetilde{N}_\gamma^\delta-\widetilde{C}_{\beta\gamma\delta}\widetilde{N}_\eta^\delta+\widetilde{C}_{\gamma\eta\delta}\widetilde{N}_\beta^\delta\right).$$
    Substituting the transformation formulas from Lemmas \ref{lem2}, \ref{lem4}, and \ref{lem8}, we obtain
    $$\widetilde\Gamma_{\beta\gamma}^\alpha=\gamma^\alpha_{\beta\gamma}-\left(g_{\beta\gamma}\sigma^\alpha-\delta^\alpha_\beta\sigma_\gamma-\delta^\alpha_\gamma\sigma_\beta\right)-g^{\alpha\eta}\left(C_{\eta\beta\delta}N_{\gamma}^{\delta}-C_{\eta\beta\delta}B_{\gamma}^{\delta}-C_{\beta\gamma\delta}N_{\eta}^{\delta}+C_{\beta\gamma\delta}B_{\eta}^{\delta}+C_{\gamma\eta\delta}N_{\beta}^{\delta}-C_{\gamma\eta\delta}B_{\beta}^{\delta}\right).$$
    The terms involving $N^\delta_\gamma$, $N^\delta_\eta$, and
    $N^\delta_\beta$ combine with $\gamma_{\beta\gamma}^\alpha$
    to recover $\Gamma_{\beta\gamma}^\alpha$ by \eqref{Gamma1}.
    The remaining terms yield
    $$\widetilde\Gamma_{\beta\gamma}^\alpha=\Gamma_{\beta\gamma}^\alpha-\left(g_{\beta\gamma}\sigma^\alpha-\delta^\alpha_\beta\sigma_\gamma-\delta^\alpha_\gamma\sigma_\beta\right)
    +g^{\alpha\eta}\left(C_{\eta\beta\delta}B_{\gamma}^{\delta}-C_{\beta\gamma\delta}B_{\eta}^{\delta}+C_{\gamma\eta\delta}B_{\beta}^{\delta}\right).$$
    This is precisely \eqref{Gamma2}, and the proof is complete.
\end{proof}

\begin{lemma}\label{lem10}
    Let $F$ and $\widetilde F=e^{\sigma}F$ be two conformally related Finsler metrics on $N$. Their Legendre transformations and inverse Legendre transformations satisfy 
    \begin{align*}
    \widetilde{\mathcal{L}}=e^{2\sigma}\mathcal{L},\quad \widetilde{\mathcal{L}}^{-1}=e^{-2\sigma}\mathcal{L}^{-1}.
    \end{align*}
\end{lemma}

\begin{proof}
    For any $(x,y)\in TN\backslash{\{0\}}$, Lemma \ref{lem1} gives
    \begin{align*}
        \widetilde{\mathcal{L}}(x,y)=\widetilde{g}_{y}(y,\cdot)=\widetilde{g}_{\alpha\beta}(x,y)y^\beta\mathrm{d}x^\alpha
        =e^{2\sigma}g_{\alpha\beta}(x,y)y^\beta\mathrm{d}x^\alpha=e^{2\sigma}\mathcal{L}(x,y).
    \end{align*}
    The same identity holds at $y=0$, since both Legendre transformations vanish there. Thus
    $\widetilde{\mathcal{L}}=e^{2\sigma}\mathcal{L}$.
    
    We now prove the corresponding identity for the inverse Legendre transformations. Let $\xi\in T_x^*N$ and set $y=\widetilde{\mathcal{L}}^{-1}(\xi)$. By the definition of the inverse Legendre transformation, we have $\xi=\widetilde{\mathcal{L}}(y)$. Using the identity established above, it follows that
    \begin{align*}
        \xi=e^{2\sigma}\mathcal{L}(y),\quad 
        \mathcal{L}(y)=e^{-2\sigma}\xi.
    \end{align*}
    Since the Legendre transformation is positively homogeneous of degree one in the fiber variable, its inverse is also positively homogeneous of degree one. Hence, applying $\mathcal{L}^{-1}$ to both sides gives
    $$y=e^{-2\sigma}\mathcal{L}^{-1}(\xi).$$
    Consequently,
    $$\widetilde{\mathcal{L}}^{-1}(\xi)=y=e^{-2\sigma}\mathcal{L}^{-1}(\xi),$$
    and thus
    $$\widetilde{\mathcal{L}}^{-1}=e^{-2\sigma}\mathcal{L}^{-1}.$$
    This completes the proof.
\end{proof}

The transformation formulas established in this section will be used in the next section to express the anisotropic mean curvature and the corresponding geometric flow in terms of the original Minkowski metric.

\section{The anisotropic mean curvature flow equation in the Gaussian space}
Throughout this section, we take $(N,F)$ to be the $(n+1)$-dimensional Minkowski space $(\mathbb{R}^{n+1},F)$. Let $$\widetilde F=e^{\sigma}F,\quad \sigma(x)=-\frac{|x|^2}{n},$$
where $|x|$ denotes the standard Euclidean norm of $x\in\mathbb{R}^{n+1}$. Then $\big(\mathbb{R}^{n+1},\widetilde F\big)$ is the Gaussian Finsler space associated with $(\mathbb{R}^{n+1},F)$.

For an immersed hypersurface $\varphi:M^n\rightarrow\mathbb{R}^{n+1}$, we denote by $H$ and $\widetilde H$ its anisotropic mean curvatures with respect to $F$ and $\widetilde F$, respectively. The corresponding unit normal vector fields are denoted by $\mathbf{n}$ and $\widetilde{\mathbf{n}}_0$.

The purpose of this section is to express the anisotropic mean curvature flow in the Gaussian space in terms of the anisotropic mean curvature and unit normal associated with the original Minkowski metric $F$. To this end, we first derive the transformation formula relating the anisotropic mean curvatures under the conformal change of the ambient Finsler metric.
 
For a fixed $F^*$-unit conormal field $\nu\in\mathcal{V}(M)$, let
$$\mathbf{n}=\mathcal{L}^{-1}(\nu),\quad \widetilde{\mathbf{n}}=\widetilde{\mathcal{L}}^{-1}(\nu).$$
By Lemma \ref{lem10}, we have
$$\widetilde{\mathbf{n}}=e^{-2\sigma}\mathbf{n}=e^{\frac{2}{n}\left|\varphi\right|^2}\mathbf{n}.$$
Here $\widetilde{\mathbf{n}}$ is the Legendre dual of the conormal $\nu$, and it is not necessarily a unit vector with respect to $\widetilde F$. The corresponding $\widetilde F$-unit normal is $$\widetilde{\mathbf{n}}_0=\frac{\widetilde{\mathbf{n}}}{\widetilde{F}(\widetilde{\mathbf{n}})}=e^\sigma\widetilde{\mathbf{n}}=e^{-\sigma}\mathbf{n}.$$
Set $\widehat{g}:=\varphi^*g_{\mathbf{n}}$ and $\widetilde{\widehat{g}}:=\varphi^*\widetilde{g}_{\widetilde{\mathbf{n}}}$. For the hypersurface $\varphi:M^n\rightarrow (\mathbb{R}^{n+1}, F)$, the Gauss formula reads 
\begin{align*}
    D^{\mathbf{n}}_{\partial_j}\varphi_*(\partial_i)=\varphi_*(\widehat\nabla_{\partial_j}\partial_i)+h_{ij}\mathbf{n},
\end{align*}
or, in local coordinates,
\begin{align}\label{Gauss1}       \partial_j\varphi^\alpha_i+\Gamma_{\beta\gamma}^\alpha(\mathbf{n})\varphi^\beta_i\varphi^\gamma_j=\widehat{\Gamma}_{ij}^k\varphi^\alpha_k+h_{ij}\mathbf{n}^\alpha. 
\end{align}
Similarly, for the hypersurface $\varphi:M^n\rightarrow (\mathbb{R}^{n+1}, \widetilde{F})$, we have 
\begin{align*}
    \widetilde{D}^{\widetilde{\mathbf{n}}}_{\partial_j}\varphi_*(\partial_i)=\varphi_*(\widetilde{\widehat\nabla}_{\partial_j}\partial_i)+\widetilde{h}_{ij}\widetilde{\mathbf{n}}_0,
\end{align*}   
which is equivalent to
\begin{align}\label{Gauss2} 
    \partial_j\varphi^\alpha_i+\widetilde{\Gamma}_{\beta\gamma}^\alpha(\widetilde{\mathbf{n}})\varphi^\beta_i\varphi^\gamma_j=\widetilde{\widehat{\Gamma}}_{ij}^k\varphi^\alpha_k+\widetilde{h}_{ij}\widetilde{\mathbf{n}}_0^{\alpha}.   
\end{align}
Since the Chern connection coefficients are homogeneous of degree zero in the reference vector, we have $\widetilde\Gamma(\widetilde{\mathbf{n}}_0)=\widetilde\Gamma(\widetilde{\mathbf{n}})=\widetilde\Gamma(\mathbf{n})$. Indeed, $\widetilde{\mathbf{n}}_0$ and $\widetilde{\mathbf{n}}$ are positively colinear, and $\widetilde{\mathbf{n}}$ and $\mathbf{n}$ are also positively colinear. Moreover, in the conformal transformation formula \eqref{Gamma2}, the Cartan tensor $C$ is homogeneous of degree $-1$, whereas $B$ is homogeneous of degree $1$. Hence the terms involving the product of $C$ and $B$ are homogeneous of degree zero, and the transformation formula can be evaluated at the reference vector $\mathbf{n}$.

By Lemma \ref{lem9}, the Chern connection coefficients transform according to
\begin{align}\label{Gamma}             
    \widetilde{\Gamma}_{\beta\gamma}^\alpha(\widetilde{\mathbf{n}})&=\Gamma_{\beta\gamma}^\alpha(\mathbf{n})-\left(g_{\beta\gamma}(\mathbf{n})g^{\alpha\kappa}(\mathbf{n})\sigma_{\kappa}-\delta_{\beta}^{\alpha}\sigma_{\gamma}-\delta_{\gamma}^{\alpha}\sigma_{\beta}\right)\nonumber\\
    &+g^{\alpha\eta}(\mathbf{n})\left(C_{\eta\beta\delta}(\mathbf{n})B_{\gamma}^{\delta}(\mathbf{n})-C_{\beta\gamma\delta}(\mathbf{n})B_{\eta}^{\delta}(\mathbf{n})+C_{\gamma\eta\delta}(\mathbf{n})B_{\beta}^{\delta}(\mathbf{n})\right),
\end{align}
where
    \begin{align*}
        B_{\beta}^{\alpha}(\mathbf{n})=g_{\beta\gamma}(\mathbf{n})\mathbf{n}^{\gamma}g^{\alpha\kappa}(\mathbf{n})\sigma_{\kappa}-\delta_{\beta}^{\alpha}\sigma_{\kappa}\mathbf{n}^{\kappa}-\mathbf{n}^{\alpha}\sigma_{\beta}-F^2g^{\alpha\lambda}(\mathbf{n})g^{\kappa\mu}(\mathbf{n})C_{\lambda\mu\beta}(\mathbf{n})\sigma_{\kappa}.
    \end{align*}
Substituting \eqref{Gamma} into the left-hand side of \eqref{Gauss2}, we obtain
\begin{align}\label{star1}
    &\partial_j\varphi^\alpha_i+\Gamma_{\beta\gamma}^\alpha(\mathbf{n})\varphi^\beta_i\varphi^\gamma_j-\left(g_{\beta\gamma}(\mathbf{n})g^{\alpha\kappa}(\mathbf{n})\sigma_{\kappa}-\delta_{\beta}^{\alpha}\sigma_{\gamma}-\delta_{\gamma}^{\alpha}\sigma_{\beta}\right)\varphi^\beta_i\varphi^\gamma_j\nonumber\\
    &+g^{\alpha\eta}(\mathbf{n})\left(C_{\eta\beta\delta}(\mathbf{n})B_{\gamma}^{\delta}(\mathbf{n})-C_{\beta\gamma\delta}(\mathbf{n})B_{\eta}^{\delta}(\mathbf{n})+C_{\gamma\eta\delta}(\mathbf{n})B_{\beta}^{\delta}(\mathbf{n})\right)\varphi^\beta_i\varphi^\gamma_j\nonumber\\
    =&\smash{\widetilde{\widehat{\Gamma}}}_{ij}^k\varphi^\alpha_k+\widetilde{h}_{ij}\widetilde{\mathbf{n}}_0^{\alpha}.     
\end{align}
We now simplify the left-hand side of \eqref{star1}. Since $\widehat{g}_{ij}=g_{\beta\gamma}(\mathbf{n})\varphi^\beta_i\varphi^\gamma_j$, we have
$$-g_{\beta\gamma}(\mathbf{n})g^{\alpha\kappa}(\mathbf{n})\sigma_{\kappa}\varphi^\beta_i\varphi^\gamma_j=-\widehat{g}_{ij}g^{\alpha\kappa}(\mathbf{n})\sigma_{\kappa}.$$
Furthermore, using $\partial_j\sigma=\sigma_\gamma\varphi^\gamma_j$, we obtain
$$\delta_{\beta}^{\alpha}\sigma_{\gamma}\varphi^\beta_i\varphi^\gamma_j=\varphi^\alpha_i\partial_j\sigma,\quad \delta_{\gamma}^{\alpha}\sigma_{\beta}\varphi^\beta_i\varphi^\gamma_j=\varphi^\alpha_j\partial_i\sigma.$$ 
Combining these identities with \eqref{Gauss1}, equation \eqref{star1} becomes
\begin{align}\label{star2}
    &\widehat{\Gamma}_{ij}^k\varphi^\alpha_k+h_{ij}\mathbf{n}^\alpha-\widehat{g}_{ij}g^{\alpha\kappa}(\mathbf{n})\sigma_{\kappa}+\varphi^\alpha_i\partial_j\sigma+\varphi^\alpha_j\partial_i\sigma\nonumber\\
    &+g^{\alpha\eta}(\mathbf{n})\left(C_{\eta\beta\delta}(\mathbf{n})B_{\gamma}^{\delta}(\mathbf{n})-C_{\beta\gamma\delta}(\mathbf{n})B_{\eta}^{\delta}(\mathbf{n})+C_{\gamma\eta\delta}(\mathbf{n})B_{\beta}^{\delta}(\mathbf{n})\right)\varphi^\beta_i\varphi^\gamma_j\nonumber\\ =&\smash{\widetilde{\widehat{\Gamma}}}_{ij}^k\varphi^\alpha_k+\widetilde{h}_{ij}\widetilde{\mathbf{n}}_0^{\alpha}.   
\end{align}

To derive the relation between the second fundamental forms $\widetilde{h}_{ij}$ and $h_{ij}$, we take the normal component of the left-hand side of \eqref{star2} with respect to $g_{\mathbf{n}}$. The first, fourth, and fifth terms are tangential, while the second term is normal. The third term admits the decomposition
\begin{align}\label{term3}
    -\widehat{g}_{ij}g^{\alpha\kappa}(\mathbf{n})\sigma_{\kappa}
    &=-\widehat{g}_{ij}(\nabla^\mathbf{n}\sigma)^\alpha
    =-\widehat{g}_{ij}\left((\nabla^\mathbf{n}\sigma)^\top\right)
    ^\alpha
    -\widehat{g}_{ij}\left((\nabla^\mathbf{n}\sigma)^\perp\right)
    ^\alpha\nonumber\\
    &=-\widehat{g}_{ij}\left((\nabla^\mathbf{n}\sigma)^\top\right)
    ^\alpha
    -\widehat{g}_{ij}\mathbf{n}(\sigma)\mathbf{n}^\alpha,
\end{align}
where $\nabla^{\mathbf{n}}\sigma=g^{\alpha\beta}(\mathbf{n})\sigma_{\beta}\frac{\partial}{\partial x^\alpha}$ is the gradient of $\sigma$ with reference vector $\mathbf{n}$.

Let $T$ denote the vector whose components are given by
\begin{align*}
    T^\alpha=g^{\alpha\eta}(\mathbf{n})\left(C_{\eta\beta\delta}(\mathbf{n})B_{\gamma}^{\delta}(\mathbf{n})-C_{\beta\gamma\delta}(\mathbf{n})B_{\eta}^{\delta}(\mathbf{n})+C_{\gamma\eta\delta}(\mathbf{n})B_{\beta}^{\delta}(\mathbf{n})\right)\varphi^\beta_i\varphi^\gamma_j.
\end{align*}
To determine its normal component, we compute
\begin{align}\label{Tn1}
    g_{\mathbf{n}}(T,\mathbf{n})&=T^\alpha g_{\alpha\kappa}(\mathbf{n})\mathbf{n}^\kappa\nonumber\\
    &=\mathbf{n}^\eta\left(C_{\eta\beta\delta}(\mathbf{n})B_{\gamma}^{\delta}(\mathbf{n})-C_{\beta\gamma\delta}(\mathbf{n})B_{\eta}^{\delta}(\mathbf{n})+C_{\gamma\eta\delta}(\mathbf{n})B_{\beta}^{\delta}(\mathbf{n})\right)\varphi^\beta_i\varphi^\gamma_j.  
\end{align}
Since $C_{\eta\beta\delta}(\mathbf{n})\mathbf{n}^\eta=0$, the first and third terms vanish. Thus 
$$g_{\mathbf{n}}(T,\mathbf{n})=-\mathbf{n}^\eta C_{\beta\gamma\delta}(\mathbf{n})B_{\eta}^{\delta}(\mathbf{n})\varphi^\beta_i\varphi^\gamma_j.$$
By the definition of $B_\eta^\delta$, we have
\begin{align*}
    \mathbf{n}^\eta B_{\eta}^{\delta}(\mathbf{n})=g^{\delta\kappa}\sigma_{\kappa}-2\mathbf{n}^\delta\sigma_0,\quad \sigma_0:=\sigma_{\kappa}\mathbf{n}^{\kappa}.  
\end{align*}
Using $C_{\beta\gamma\delta}(\mathbf{n})\mathbf{n}^\delta=0$, it follows that
\begin{align*}
    g_{\mathbf{n}}(T,\mathbf{n})=-C_{\beta\gamma\delta}(\mathbf{n})g^{\delta\kappa}\sigma_{\kappa}\varphi^\beta_i\varphi^\gamma_j=-C_{\beta\gamma}^\kappa(\mathbf{n})\sigma_{\kappa}\varphi^\beta_i\varphi^\gamma_j.
\end{align*}
Consequently, the normal component of $T$ is
\begin{align}\label{Tn2}
    (T^\perp)^\alpha=g_{\mathbf{n}}(T,\mathbf{n})\mathbf{n}^\alpha=-C_{\beta\gamma}^\kappa(\mathbf{n})\sigma_{\kappa}\varphi^\beta_i\varphi^\gamma_j\mathbf{n}^\alpha.
\end{align}
Combining \eqref{term3} and \eqref{Tn2}, the normal component of the left-hand side of \eqref{star2} is
\begin{align}\label{bot1}
    &h_{ij}\mathbf{n}^\alpha-\widehat{g}_{ij}\left(g^{\alpha\kappa}(\mathbf{n})\sigma_\kappa\right)^\bot+(T^\bot)^\alpha\nonumber\\
    =&\left(h_{ij}-\widehat{g}_{ij}\mathbf{n}(\sigma)-C_{\beta\gamma}^\kappa(\mathbf{n})\sigma_{\kappa}\varphi^\beta_i\varphi^\gamma_j\right)\mathbf{n}^\alpha,
\end{align}
On the other hand, the normal component of the right-hand side of \eqref{star2} is
\begin{align}\label{bot2}
\widetilde{h}_{ij}\widetilde{\mathbf{n}}_0^\alpha=\widetilde{h}_{ij}e^{-\sigma}\mathbf{n}^\alpha.
\end{align}
Comparing \eqref{bot1} and \eqref{bot2}, we obtain
\begin{align}
    \widetilde{h}_{ij}=e^{\sigma}\left(h_{ij}-\widehat{g}_{ij}\mathbf{n}(\sigma)-C_{\beta\gamma}^\kappa(\mathbf{n})\sigma_{\kappa}\varphi^\beta_i\varphi^\gamma_j\right).
\end{align}
Since $\widetilde{\widehat{g}}=e^{2\sigma}\widehat{g}$, taking the trace with respect to $\widetilde{\widehat{g}}$ yields
\begin{align}\label{H1}
\widetilde{H}=\smash{\widetilde{\widehat{g}}}^{ij}\widetilde{h}_{ij}=e^{-\sigma}\left(H-n\mathbf{n}(\sigma)-\widehat{I}(\nabla^\mathbf{n}\sigma)\right),
\end{align}
where $\widehat{I}(\nabla^\mathbf{n}\sigma)=\widehat{g}^{ij}C_{\mathbf{n}}(\varphi_*(\partial_i),\varphi_*(\partial_j),\nabla^\mathbf{n}\sigma)$
is the mean Cartan torsion evaluated at $\nabla^\mathbf{n}\sigma$.
Substituting $\sigma=-\frac{|\varphi|^2}{n}$ into \eqref{H1}, we obtain
\begin{align}\label{H2}
\widetilde{H}=e^{\frac{|\varphi|^2}{n}}\left(H+2\langle\varphi,\mathbf{n}\rangle+\frac{1}{n}\widehat{I}(\nabla^\mathbf{n}|\varphi|^2)\right),
\end{align}
where $\widehat{I}(\nabla^\mathbf{n}|\varphi|^2)
=\widehat{g}^{ij}C_{\mathbf{n}}(\varphi_*(\partial_i),\varphi_*(\partial_j),\nabla^\mathbf{n}|\varphi|^2)
=2\widehat{g}^{ij}C_{\beta\gamma}^{\kappa}(\mathbf{n})\varphi^\kappa\varphi^\beta_i\varphi^\gamma_j$.

We are now in a position to express the anisotropic mean curvature flow in the Gaussian space in terms of the anisotropic mean curvature $H$ and the unit normal $\mathbf{n}$ associated with the original Minkowski metric $F$. This leads to the following definition.
\begin{definition}
    Let $\varphi_0: M^n\rightarrow (\mathbb{R}^{n+1},\widetilde F)$ be a smooth immersion of an $n$-dimensional manifold $M^n$, and let $\widetilde{\widehat{g}}_0:=\varphi_0^*\widetilde{g}_{\widetilde{\mathbf{n}}_0}$ be the induced Riemannian metric on $M^n$. The anisotropic mean curvature flow of $\varphi_0$ in the Gaussian space is a family of smooth immersions $\varphi_t: M\rightarrow (\mathbb{R}^{n+1},\widetilde F)$, $t\in[0,T)$, such that $\varphi:M\times[0,T)\to(\mathbb{R}^{n+1},\widetilde F)$ satisfies
    \begin{align}\label{mainflow}
        \begin{cases}
            \partial_t\varphi(p,t)=\widetilde H(p,t)\widetilde{\mathbf{n}}_0(p,t)
            =e^{\frac{2}{n}\left|\varphi\right|^2}\left(H+2\langle\varphi,\mathbf{n}\rangle+\frac{1}{n}\widehat{I}(\nabla^\mathbf{n}|\varphi|^2)\right)\mathbf{n},\\
            \varphi(p,0)=\varphi_0(p),
        \end{cases}
    \end{align}
    where $H$ and $\mathbf{n}$ denote, respectively, the anisotropic mean curvature and the inward-pointing unit normal of $\varphi_t$ with respect to $F$.
\end{definition}

\begin{proposition}\label{prop+1}
    Let $\varphi_t:(M^n,\widehat g)\to(\mathbb{R}^{n+1},F)$ be a smooth family of immersions and $V(p,t)\in \mathrm{d}\varphi_t(T_pM)$ be a smooth vector field along $M^n$ for all $(p,t)\in M^n\times[0,T)$. If $\varphi_t$ satisfies
    \begin{align}\label{+1}
        \begin{cases}
            \partial_t\varphi(p,t)=e^{\frac{2}{n}\left|\varphi(p,t)\right|^2}\left(H(p,t)+2\langle\varphi(p,t),\mathbf{n}(p,t)\rangle+\frac{1}{n}\widehat{I}(\nabla^\mathbf{n}|\varphi(p,t)|^2)\right)\mathbf{n}(p,t)+V(p,t),\\
            \varphi(p,0)=\varphi_0(p),
        \end{cases}
    \end{align}
    then there locally exists a family of reparametrizations of $\varphi_t$ which satisfies \eqref{mainflow} around any point of $M^n\times[0,T)$.

    In particular, if $M^n$ is compact, then there uniquely exists a family of global reparametrizations of $\varphi_t$ as above for every $t\geq0$ which satisfies \eqref{mainflow} and leaves $\varphi_0$ unmodified.

    On the other hand, if $\varphi_t:(M^n,\widehat g)\to(\mathbb{R}^{n+1},F)$ can be globally reparametrized for $t\geq0$ such that it satisfies \eqref{mainflow}, then $\varphi$ must satisfy \eqref{+1} for some vector field $V(q,t)\in\mathrm{d}\varphi_t(T_pM^n)$.
\end{proposition}

\begin{proof}
    First, we assume that $M^n$ is compact. The proof here follows almost [\citenum{mc11}, Proposition 1.3.4] and [\citenum{zhc18}, Proposition 3.2]. We write it for the reader's convenience and dispelling possible doubts.

    Let $\check{\varphi}(p,t)=\varphi(\Phi(p,t),t)$ and $\Phi:M^n\times[0,T)\to M^n$ be a smooth family of diffeomorphisms of $M^n$ with $\Phi(p,0)=p$ for every $p\in M^n$ and 
    \begin{align}\label{+2}
        \partial_t\Phi(p,t)=W(\Phi(p,t),t),
    \end{align}
    where $W(q,t)=-[\mathrm{d}\varphi_t]^{-1}(V(q,t))$ is a time-dependent vector field on $M^n$. By the existence and uniqueness theorem for ODEs on the compact manifold $M^n$, this family uniquely and smoothly exists.

    Considering the reparametrizations $\check{\varphi}(p,t)$, we have
    \begin{align*}
        \check{\varphi}_i(p,t)
        &=\varphi_k(\Phi(p,t),t)\Phi^k_i,\\
        \check{\varphi}_{ij}(p,t)
        &=\varphi_{kl}(\Phi(p,t),t)\Phi^k_i\Phi^l_j
        +\varphi_k(\Phi(p,t),t)\Phi^k_{ij},\\
        \check{\mathbf{n}}(p,t)
        &=\mathbf{n}(\Phi(p,t),t),\\
        \check{g}_{ij}(p,t)
        &=\widehat{g}_{kl}(\Phi(p,t),t)\Phi^k_i\Phi^l_j,\\
        \check{g}^{ij}(p,t)&\Phi^k_i\Phi^l_j
        =\widehat{g}^{kl}(\Phi(p,t),t),\\
        \check{h}_{ij}(p,t)
        &=g_{\check{\mathbf{n}}}\big(\check{\mathbf{n}},(\check{\varphi}^{\alpha}_{ij}+\Gamma_{\beta\gamma}^{\alpha}(\check{\mathbf{n}})\check{\varphi}^{\beta}_i\smash{\check{\varphi}}^{\gamma}_j)\partial_{\alpha}\big)\\
        &=g_{\mathbf{n}}\big(\mathbf{n},(\varphi^{\alpha}_{kl}+\Gamma_{\beta\gamma}^{\alpha}(\mathbf{n})\varphi^{\beta}_k\varphi^{\gamma}_l)\partial_{\alpha}\big)\\
        &=h_{kl}(\Phi(p,t),t)\Phi^k_i\Phi^l_j,\\
        \check{H}(p,t)&=H(\Phi(p,t),t),
    \end{align*}
    from which one has
    \begin{align*}
        \partial_t\check{\varphi}(p,t)=&\partial_t\varphi(\Phi(p,t),t)+\mathrm{d}\varphi(\Phi(p,t))(\partial_t\Phi(p,t))\\
        =&e^{\frac{2}{n}\left|\varphi(\Phi(p,t),t)\right|^2}\Big(H(\Phi(p,t),t)+2\langle\varphi(\Phi(p,t),t),\mathbf{n}(\Phi(p,t),t)\rangle\\
        &+\frac{1}{n}\widehat{I}(\nabla^\mathbf{n}|\varphi(\Phi(p,t),t)|^2)\Big)\mathbf{n}(\Phi(p,t),t)+V(\Phi(p,t),t)+\mathrm{d}\varphi(\Phi(p,t))(W(\Phi(p,t),t))\\
        =&e^{\frac{2}{n}\left|\varphi(\Phi(p,t),t)\right|^2}\Big(H(\Phi(p,t),t)+2\langle\varphi(\Phi(p,t),t),\mathbf{n}(\Phi(p,t),t)\rangle\\
        &+\frac{1}{n}\widehat{I}(\nabla^\mathbf{n}|\varphi(\Phi(p,t),t)|^2)\Big)\mathbf{n}(\Phi(p,t),t)
        +V(\Phi(p,t),t)-\mathrm{d}\varphi(\Phi(p,t))[\mathrm{d}\varphi(\Phi(p,t))]^{-1}(V(\Phi(p,t),t))\\
        =&e^{\frac{2}{n}\left|\check{\varphi}(p,t)\right|^2}\Big(\check{H}(p,t)+2\langle\check{\varphi}(p,t),\check{\mathbf{n}}(p,t)\rangle
        +\frac{1}{n}\check{I}(\nabla^{\check{\mathbf{n}}}|\check{\varphi}(p,t)|^2)\Big)\check{\mathbf{n}}(p,t),
    \end{align*}
    where we have used that
    \begin{align*}
        \check{I}(\nabla^{\check{\mathbf{n}}}|\check{\varphi}(p,t)|^2)=&\check{g}^{ij}C_{\check{\mathbf{n}}}(\check{\varphi}_*(\partial_i),\check{\varphi}_*(\partial_j),\nabla^{\check{\mathbf{n}}}|\check{\varphi}|^2)|_{(p,t)}=2\check{g}^{ij}C_{\beta\gamma}^{\kappa}(\check{\mathbf{n}})\check{\varphi}^{\kappa}\check{\varphi}^{\beta}_i\check{\varphi}^{\gamma}_j|_{(p,t)}\\
        =&2\check{g}^{ij}(p,t)C_{\beta\gamma}^{\kappa}(\mathbf{n}(\Phi(p,t),t))\varphi^{\kappa}(\Phi(p,t),t)\varphi^{\beta}_k(\Phi(p,t),t)\Phi^k_i\varphi^{\gamma}_l(\Phi(p,t),t)\Phi^l_j\\
        =&2\widehat{g}^{ij}(\Phi(p,t),t)C_{\beta\gamma}^{\kappa}(\mathbf{n}(\Phi(p,t),t))\varphi^{\kappa}(\Phi(p,t),t)\varphi^{\beta}_k(\Phi(p,t),t)\varphi^{\gamma}_l(\Phi(p,t),t)\\
        =&\widehat{I}(\nabla^{\mathbf{n}}|\varphi(\Phi(p,t),t)|^2).
    \end{align*}
    Hence, $\check{\varphi}$ satisfies system \eqref{mainflow}, and $\check{\varphi}_0=\varphi_0$.

    In turn, the above computations also show that if $\check{\varphi}(p,t)=\varphi(\Phi(p,t),t)$ satisfies \eqref{mainflow}, $\Phi:M^n\times[0,T)\to M^n$ must solve \eqref{+2}. Hence, if we assume $\Phi(\cdot,0)=\mathrm{Id}_{M^n}$ in order that $\varphi_0$ is unmodified, it is unique.

    In the noncompact case, we have to solve \eqref{+2} in $[t_1,t_2]\times\Sigma$, where $0<t_1<t_2<T$ and $\Sigma\subset M^n$ with compact closure, and then we obtain a solution of \eqref{mainflow} in a possibly smaller open subset of $\Sigma$ and some interval of time.

    Now assume that $\check{\varphi}(p,t)=\varphi(\Phi(p,t),t)$ satisfies \eqref{mainflow}. Differentiating it, we have
    \begin{align*}
        \partial_t\check{\varphi}(p,t)=&\partial_t\varphi(\Phi(p,t),t)+\mathrm{d}\varphi(\Phi(p,t))(\partial_t\Phi(p,t))\\
        =&e^{\frac{2}{n}\left|\check{\varphi}(p,t)\right|^2}\Big(\check{H}(p,t)+2\langle\check{\varphi}(p,t),\check{\mathbf{n}}(p,t)\rangle
        +\frac{1}{n}\check{I}(\nabla^{\check{\mathbf{n}}}|\check{\varphi}(p,t)|^2)\Big)\check{\mathbf{n}}(p,t)\\
        =&e^{\frac{2}{n}\left|\varphi(\Phi(p,t),t)\right|^2}\Big(H(\Phi(p,t),t)+2\langle\varphi(\Phi(p,t),t),\mathbf{n}(\Phi(p,t),t)\rangle\\
        &+\frac{1}{n}\widehat{I}(\nabla^\mathbf{n}|\varphi(\Phi(p,t),t)|^2)\Big)\mathbf{n}(\Phi(p,t),t),
    \end{align*}
    i.e.,
    \begin{align*}
        \partial_t\varphi(q,t)=&e^{\frac{2}{n}\left|\varphi(q,t)\right|^2}\Big(H(q,t)+2\langle\varphi(q,t),\mathbf{n}(q,t)\rangle
        +\frac{1}{n}\widehat{I}(\nabla^\mathbf{n}|\varphi(q,t)|^2)\Big)\mathbf{n}(q,t)\\
        &-\mathrm{d}\varphi(q)
        \big(\partial_t\Phi(\Phi_t^{-1}(q),t)\big),
    \end{align*}
    for every $q\in M^n$ and $t\in[0,T)$. Then, the last statement of the proposition follows by setting $V(q,t)=-\mathrm{d}\varphi(q)
        \big(\partial_t\Phi(\Phi_t^{-1}(q),t)\big)$.
\end{proof}

\begin{corollary}\label{coro+1}
    If $\varphi_t=\varphi(\cdot,t)$ satisfies $g_{\mathbf{n}}(\partial_t\varphi,\mathbf{n})=e^{\frac{2}{n}\left|\varphi\right|^2}\Big(H+2\langle\varphi,\mathbf{n}\rangle+\frac{1}{n}\widehat{I}(\nabla^\mathbf{n}|\varphi|^2)\Big)$, then there is a (time-dependent) reparametrization $\check{\varphi}_t$ of $M^n$, such that $\varphi_t\circ\check{\varphi}_t$ satisfies \eqref{mainflow}. If $M^n$ is compact, this can be done uniquely by global reparametrizations, without modifying $\varphi_0$.
\end{corollary}

\section{Short-time existence of the flow}
In this section, we establish the short-time existence and uniqueness of the anisotropic mean curvature flow \eqref{mainflow} in the Gaussian space.
\begin{theorem}\label{mainthm}
The system \eqref{mainflow} has a unique, smooth solution $\varphi_t(M^n)$ for a short time with any smooth compact initial hypersurface $\varphi_0(M^n)$ of $\mathbb{R}^{n+1}$ at $t=0$.
\end{theorem}

\begin{proof}
    Let $\varphi_{0}:M^n\rightarrow\mathbb{R}^{n+1}$ be a smooth immersion of a compact manifold. Since the nonembedded case is similar, we now only assume that this hypersurface is embedded, which implies the inner pointing unit normal vector field $\overline{\mathbf{n}}_0$ (with respect to the Euclidean metric) is globally defined and smooth.

    We look for a smooth solution $\varphi:M^n\times[0,T)\rightarrow\mathbb{R}^{n+1}$ of the parabolic problem \eqref{mainflow} for some $T>0$. Since we are interested in a solution for short time, we can forget about the immersion condition ($d\varphi_{t}$ nonsingular) as it will follow automatically by the smoothness of the solution and by the fact that $\varphi_{0}$ is a compact immersion, when $t$ is close to zero.

    By Proposition \ref{prop+1} and Corollary \ref{coro+1}, it suffices to find a smooth solution of
    \begin{align}\label{eq1}
        \begin{cases}
            g_{\mathbf{n}}\left(\partial_{t}\varphi(p,t),\mathbf{n}(p,t)\right)=e^{\frac{2}{n}\left|\varphi\right|^2}\left(H+2\langle\varphi,\mathbf{n}\rangle+\frac{1}{n}\widehat I(\nabla^\mathbf{n}|\varphi|^2)\right)(p,t),\\
            \varphi(p,0)=\varphi_0(p).
        \end{cases}
    \end{align}

    We consider the regular Euclidean tubular neighborhood $\Gamma=\left\{x\in \mathbb{R}^{n+1}|d\left(x,\varphi_{0}(M^{n})\right)<\varepsilon\right\}$, which exists for $\varepsilon>0$ small enough. By regular we mean that the map $\Upsilon:M^n\times(-\varepsilon,\varepsilon)\rightarrow\Gamma$ defined as $\Upsilon(p,s):=\varphi_0(p)+s\overline{\mathbf{n}}_0(p)$ is a diffeomorphism, where $\overline{\mathbf{n}}_0$ is the Euclidean unit normal vector field along $\varphi_0$. Any sufficiently small deformation of $\varphi_0(M^n)$ inside $\Gamma$ can be represented as the graph of a height function $f$ over $\varphi_0(M^n)$. Conversely, to any sufficiently small function $f:M^n\rightarrow\mathbb{R}$, we can associate the hypersurface $M^f\subset\Gamma$ given by $\varphi(p)=\varphi_0(p)+f(p)\overline{\mathbf{n}}_0(p)$. We want to compute now the equation for a smooth time-dependent function $f$, such that $\varphi$ satisfies system \eqref{eq1}. Obviously, since $f(\cdot,0)$ gives the initial hypersurface $\varphi_0$, we have $f(p,0)=0$ for every $p\in M^n$.

    First, we compute the metric and the normal of the perturbed hypersurfaces. In local coordinates on $M$, write $\varphi_i=\partial_i\varphi$, $f_i=\partial_if$. Using the Euclidean Weingarten formula $\overline{\mathbf{n}}_{0i}(p)=-\smash{\overline{h}}_i^k(p,0)\varphi_{0k}(p,0)$, we obtain
    \begin{align}\label{varphi_i}
        \varphi_i(p,t)=\varphi_{0i}(p)+\overline{\mathbf{n}}_0(p)f_i(p,t)-f(p,t)\smash{\overline{h}}_i^k(p,0)\varphi_{0k}(p,0):=\varphi_i(p,f,\partial{f})
    \end{align}
    where $\partial{f}=(f_1,\cdots,f_n)$. The induced metric of $\varphi$ with respect to the Euclidean metric $\langle\cdot,\cdot\rangle$ in $\mathbb{R}^{n+1}$ is therefore
    \begin{align*}
        \overline{g}_{ij}(p,t)&=\langle\varphi_i(p,t),\varphi_j(p,t)\rangle\\
        &=\overline{g}_{ij}(p,0)-2f(p,t)\overline{h}_{ij}(p,0)
        +f^2(p,t)\overline{h}_{ik}\overline{g}^{kl}\overline{h}_{lj}(p,0)+f_i(p,t)f_j(p,t)\\
        &:=\overline{g}_{ij}(p,f,\partial{f}).
    \end{align*}
    In what follows, $\overline{g}^{ij}(p,t)$ denotes the inverse matrix of $\overline{g}_{ij}(p,t)$, and we denote by $\overline\nabla f=f_i\overline{g}^{ij}(p,t)\partial_j$ the gradient of $f$ with respect to the time-dependent metric $\overline g(t)$. The Euclidean normal vector of $\varphi$, which need not be normalized, can then be written as
    \begin{align*}
        \overline{\mathbf{n}}&=\overline{\mathbf{n}}_0-\mathrm d\varphi(\overline\nabla{f}).
    \end{align*}
    Indeed, this vector is orthogonal to every $\varphi_i$, since $\langle\overline{\mathbf{n}},\varphi_i\rangle=f_i-f_j\overline{g}^{jk}(t)\overline{g}_{ki}(t)=0$.
    Hence, the $F^*$-unit conormal and the corresponding $F$-unit normal are given by
    $$\nu=\frac{\overline{\mathbf{n}}}{F^*(\overline{\mathbf{n}})}=\frac{\overline{\mathbf{n}}_0-d\varphi(\overline\nabla{f})}{F^*\big(\overline{\mathbf{n}}_0-d\varphi(\overline\nabla{f})\big)}:=\nu(p,f,\partial{f}),$$
    and
    $$\mathbf{n}=\mathcal{L}^{-1}(\nu):=\mathbf{n}(p,f,\partial{f}).$$
    Consequently, the induced Finsler metric
    $$\widehat{g}_{ij}(p,t):=g_{\mathbf{n}}(\varphi_i(p,t),\varphi_j(p,t)):=\widehat{g}_{ij}(p,f,\partial{f}),$$
    and its inverse also depend smoothly on $p$, $f$, and $\partial f$, provided that the graph remains sufficiently small.
    
    Notice that the Euclidean normal, the Finsler normal, and the induced metrics, as well as their inverses, depend only on the first space derivatives of $f$. Moreover, since $f(p,0)=\overline\nabla{f}(p,0)=0$, all these quantities are smooth near the initial data. Since $M^n$ is compact, when $t$ is small the denominator $F^*(\overline{\mathbf{n}})$ is uniformly bounded below away from zero.

    From \eqref{varphi_i}, we have
    \begin{align*}
        \varphi_{ij}=\varphi_{0ij}+\overline{\mathbf{n}}_0f_{ij}+f\overline{\mathbf{n}}_{0ij}
        +f_i\overline{\mathbf{n}}_{0j}+f_j\overline{\mathbf{n}}_{0i}.
    \end{align*}
    Then, substituting $\varphi_{ij}$ into the expression for $h_{ij}$, and using $\langle\varphi_i,\overline{\mathbf{n}}_0\rangle=f_i$ and $\langle\mathrm{d}\varphi(\overline\nabla f),\overline{\mathbf{n}}_0\rangle=|\overline\nabla f|^2_{\overline g}$, we obtain
    \begin{align*}
        h_{ij}=&g_{\mathbf{n}}(\mathbf{n},\varphi_{ij})=\nu(\varphi_{ij})
        =\frac{1}{F^*({\overline{\mathbf{n}})}}\overline{\mathbf{n}}(\varphi_{ij})\\
        =&\frac{1}{F^*\big(\overline{\mathbf{n}}_0-\mathrm{d}\varphi(\overline\nabla{f})\big)}\langle\overline{\mathbf{n}}_0-\mathrm d\varphi(\overline\nabla{f}),\varphi_{ij}\rangle\\
        =&\frac{1}{F^*\big(\overline{\mathbf{n}}_0-d\varphi(\overline\nabla{f})\big)}\Big(\overline{h}_{ij}(p,0)+f_{ij}+f\langle\overline{\mathbf{n}}_{0},\overline{\mathbf{n}}_{0ij}\rangle-\langle\mathrm {d}\varphi(\overline\nabla{f}),\varphi_{0ij}\rangle\\
        &-f_{ij}\langle\mathrm{d}\varphi(\overline\nabla f),\overline{\mathbf{n}}_0\rangle-f_i\langle\mathrm{d}\varphi(\overline\nabla{f}),\overline{\mathbf{n}}_{0j}\rangle-f_j\langle\mathrm{d}\varphi(\overline\nabla{f}),\overline{\mathbf{n}}_{0i}\rangle-\langle\mathrm{d}\varphi(\overline\nabla{f}),f\overline{\mathbf{n}}_{0ij}\rangle\Big)\\
        =&\frac{1-|\overline\nabla f|^2_{\overline g}}{F^*\big(\overline{\mathbf{n}}_0-\mathrm d\varphi(\overline\nabla{f})\big)}f_{ij}+P_{ij}(p,f,\partial{f}),
    \end{align*}
    where $P_{ij}$ is a smooth form when $f$ and $\partial{f}$ are small, hence for $t$ small.

    The mean curvature is then given by 
    \begin{align*}
        H=\widehat{g}^{ij}h_{ij}=\frac{1-|\overline\nabla f|^2_{\overline g}}{F^*\big(\overline{\mathbf{n}}_0-d\varphi(\overline\nabla{f})\big)}\widehat{g}^{ij}f_{ij}+\mathcal{P}(p,f,\partial{f}),
    \end{align*}
    where $(\widehat{g}^{ij})=(\widehat{g}_{ij})^{-1}=(\widehat{g}^{ij}(p,f,\partial{f}))$, and $\mathcal{P}$ is a smooth function, assuming that $f$ and $\overline\nabla{f}$ are small.

    We are now ready to write down the condition 
    $$g_{\mathbf{n}}\left(\partial_{t}\varphi,\mathbf{n}\right)=e^{\frac{2}{n}\left|\varphi\right|^2}\Big(H+2\langle\varphi,\mathbf{n}\rangle+\frac{1}{n}\widehat I(\nabla^\mathbf{n}|\varphi|^2)
    \Big)$$
    in terms of the function $f$,
    \begin{align}\label{f}
        (\partial_{t}f)g_{\mathbf{n}}(\mathbf{n},\overline{\mathbf{n}}_0)
        =&g_{\mathbf{n}}\left(\partial_{t}\varphi(p,t),\mathbf{n}(p,t)\right)\nonumber\\
        =&e^{\frac{2}{n}\left|\varphi\right|^2}
        \left(H+2\langle\varphi,\mathbf{n}\rangle+\frac{1}{n}\widehat I(\nabla^\mathbf{n}|\varphi|^2)\right)(p,t)\nonumber\\
        =&e^{\frac{2}{n}\left|\varphi\right|^2}
        \left(H+2\langle\varphi,\mathbf{n}\rangle+\frac{2}{n}\widehat{g}^{ij}C_{\beta\gamma}^{\kappa}(\mathbf{n})\varphi^{\kappa}\varphi^{\beta}_{i}\varphi^{\gamma}_{j}\right)\nonumber\\
        =&e^{\frac{2}{n}\left|\varphi_0+f\overline{\mathbf{n}}_0\right|^2}\frac{1-|\overline\nabla f|^2_{\bar g}}{F^*\big(\overline{\mathbf{n}}_0-d\varphi(\overline\nabla{f})\big)}\widehat{g}^{ij}f_{ij}+\mathcal{R}(p,f,\partial{f}),
    \end{align}
    where 
    $$\mathcal{R}(p,f,\partial{f})=e^{\frac{2}{n}\left|\varphi_0+f\overline{\mathbf{n}}_0\right|^2}\left(\mathcal{P}(p,f,\partial{f})+2\langle\varphi,\mathbf{n}\rangle+\frac{2}{n}\widehat{g}^{ij}C_{\beta\gamma}^{\kappa}(\mathbf{n})\varphi^{\kappa}\varphi^{\beta}_{i}\varphi^{\gamma}_{j}\right).$$
    Notice that 
    $$g_{\mathbf{n}}(\mathbf{n},\overline{\mathbf{n}}_0)=\nu(\overline{\mathbf{n}}_0)=\frac{1-|\overline\nabla f|^2_{\overline g}}{F^*\big(\overline{\mathbf{n}}_0-d\varphi(\overline\nabla{f})\big)}.$$ 
    At the initial time $t=0$, where $f=0$ and $\overline\nabla f=0$, we have $g_{\mathbf{n}}(\mathbf{n},\overline{\mathbf{n}}_0)=\nu(\overline{\mathbf{n}}_0)=\frac{1}{F^*(\overline{\mathbf{n}}_0)}>0$. Since $M$ is compact, this quantity remains uniformly positive for sufficiently small $f$ and $\overline\nabla f$. Thus, dividing both sides of \eqref{f} by $g_{\mathbf{n}}(\mathbf{n},\overline{\mathbf{n}}_0)=\nu(\overline{\mathbf{n}}_0)$ is valid. We obtain
    \begin{align*}
        \partial_{t}f=&e^{\frac{2}{n}\left|\varphi_0+f\overline{\mathbf{n}}_0\right|^2}\widehat{g}^{ij}f_{ij}
        +\frac{\mathcal{R}(p,f,\partial{f})}{\nu(\overline{\mathbf{n}}_0)}
        =e^{\frac{2}{n}\left|\varphi_0+f\overline{\mathbf{n}}_0\right|^2}\widehat{g}^{ij}f_{ij}
        +\mathcal{Q}(p,f,\partial{f}),
    \end{align*}
    where $\mathcal{Q}(p,\cdot,\cdot)$ is a smooth function when its arguments are small. Moreover, we verify the strict parabolicity of this equation. At the initial time $t=0$, where $f=0$ and $\overline\nabla f=0$, the principal coefficient matrix is
    $$a^{ij}(p,0,0)=e^{\frac{2}{n}|\varphi_0|^2}\widehat{g}^{ij}(p,0).$$
    Since $\widehat{g}^{ij}(p,0)$ is the inverse of the induced Finsler metric on the initial hypersurface, it is positive definite. Since $M$ is compact, $(a^{ij}(p,0,0))$ is uniformly positive definite. By the smooth dependence of the coefficients on $p$, $f$ and $\partial f$, the principal coefficient remains uniformly positive definite whenever $f$ and $\partial f$ are sufficiently small. Hence the evolution equation is a uniformly strictly parabolic quasilinear equation in a neighborhood of the initial data.

    Therefore, if the smooth function $f:M^n\times [0,T)\rightarrow\mathbb{R}$ solves the following partial differential equation
    \begin{align}\label{ptf}
        \begin{cases}
            \partial_tf=e^{\frac{2}{n}\left|\varphi_0+f\overline{\mathbf{n}}_0\right|^2}\widehat{g}^{ij}f_{ij}
            +\mathcal{Q}(p,f,\partial{f}),\\
            f(p,0)=0,
        \end{cases}
    \end{align}
    then $\varphi(p,t)=\varphi_0(p)+f(p,t)\overline{\mathbf{n}}_0(p)$ is a solution of \eqref{mainflow} for the initial compact embedding $\varphi_0$, in a positive time interval.

    Conversely, if $\varphi$ satisfies \eqref{mainflow}, for small time the hypersurfaces $\varphi_t$ are embedded in the tubular neighborhood $\Gamma$ of $\varphi_0(M^n)$, then the function $f(p,t)=\pi_{(-\varepsilon,\varepsilon)}\left[\Upsilon^{-1}(\varphi(p,t))\right]$ is smooth and $f(p,0)=\pi_{(-\varepsilon,\varepsilon)}\left[\Upsilon^{-1}(\varphi(p,0))\right]=\pi_{(-\varepsilon,\varepsilon)}\left[(p,0)\right]=0$,
    where $\pi_{(-\varepsilon,\varepsilon)}$ is the projection map on the second factor of $M^n\times(-\varepsilon,\varepsilon)$, hence, by the above computations $f$ must solve problem \eqref{ptf}. 

    Using arguments similar to those in [\citenum{mc11}, Theorem 1.5.1, p.22], with a minor modification, we obtain the existence, uniqueness and continuous dependence of a solution for \eqref{ptf}. By the unique solution $f$ of \eqref{ptf}, we can prove existence, uniqueness and continuous dependence of a solution for \eqref{eq1}.
\end{proof}

\section*{Acknowledgments}
The author would like to express sincere gratitude to Professor Q. L. Xia for her valuable suggestions and constant encouragement. This work was supported by Zhejiang Provincial Natural Science Foundation of China (Grant No. LQ24A010020).

\end{document}